\documentclass[12pt]{article}
\usepackage{amsmath,amssymb,amsthm}
\RequirePackage[dvips]{graphicx}
\usepackage{xcolor}
\usepackage{tikz}
\usepackage{cite}

\theoremstyle{plain}
\newtheorem{thm}{Theorem}[section]

\newtheorem{lem}[thm]{Lemma}
\newtheorem{cor}[thm]{Corollary}
\newtheorem{prop}[thm]{Proposition}

\theoremstyle{definition}
\newtheorem{rem}[thm]{Remark}

\newtheorem{exam}[thm]{Example}

\DeclareMathOperator{\diag}{diag}
\DeclareMathOperator{\spec}{Spec}
\DeclareMathOperator{\lspec}{L-spec}
\DeclareMathOperator{\qspec}{Q-spec}
\DeclareMathOperator{\cnspec}{CN-spec}
\DeclareMathOperator{\cnlspec}{CNL-spec}
\DeclareMathOperator{\cnqspec}{CNSL-spec}
\DeclareMathOperator{\CNRS}{CNRS}
\DeclareMathOperator{\CNL}{CNL}
\DeclareMathOperator{\CNSL}{CNSL}
\DeclareMathOperator{\CN}{CN}

\begin{document}

  \begin{center}
{\Large \bf Common neighborhood energies and their relations with Zagreb index}\\
 \vspace{6mm}

 \bf {Firdous Ee Jannat$^a$, Rajat  Kanti Nath$^a$,  Kinkar Chandra Das$^{b,}$\footnote{Corresponding author}}

 \vspace{5mm}

{\it $^a$Department of Mathematical Science, Tezpur  University,\\
 Napaam -784028, Sonitpur, Assam, India\/}\\
 {\tt E-mail: firdusej@gmail.com,  rajatkantinath@yahoo.com} \\[2mm]

{\it $^b$Department of Mathematics, Sungkyunkwan University, \\
 Suwon 16419, Republic of Korea\/} \\
{\tt E-mail: kinkardas2003@gmail.com}

 \hspace*{30mm}

 \end{center}
 \baselineskip=0.23in

\begin{abstract}
In this paper we establish connections between   common neighborhood Laplacian and common neighborhood signless Laplacian  energies   and the first Zagreb index of a graph $\mathcal{G}$.
We  introduce the  concepts of CNL-hyperenergetic and CNSL-hyperenergetic graphs and showed that $\mathcal{G}$   is neither CNL-hyperenergetic nor CNSL-hyperenergetic if $\mathcal{G}$ is a complete bipartite graph.
We obtain certain relations between  various energies of a graph.  Finally, we conclude the paper with several bounds for common neighborhood Laplacian and signless Laplacian energies of a graph.\\

\end{abstract}

\baselineskip=0.30in

\section{Introduction} \label{S:intro}

Let $\mathcal{G}$ be a finite simple graph with adjacency matrix $A(\mathcal{G})$ and degree matrix $D(\mathcal{G})$. The spectrum of $\mathcal{G}$ is the set of all the eigenvalues of $A(\mathcal{G})$ with multiplicity.
 Let $L(\mathcal{G}) := D(\mathcal{G}) - A(\mathcal{G})$ and $Q(\mathcal{G}) := D(\mathcal{G}) + A(\mathcal{G})$ be the Laplacian and signless Laplacian matrices of $\mathcal{G}$, respectively. More than thousands of research papers are studied on the adjacency, Laplacian and signless Laplacian eigenvalues of graphs, especially the most
recent works \cite{CH,LZD,RDM} and the references therein. Then the Laplacian spectrum and the signless Laplacian spectrum of $\mathcal{G}$ are the sets of all the eigenvalues (with multiplicity) of $L(\mathcal{G})$ and $Q(\mathcal{G})$, respectively. Corresponding to these spectra, energy ($E(\mathcal{G})$), Laplacian energy ($LE(\mathcal{G})$) and signless Laplacian energy ($LE^+(\mathcal{G})$) of $\mathcal{G}$ are defined as
\[
E(\mathcal{G}) = \sum_{\alpha \in \spec(\mathcal{G})} |\alpha|, \quad LE(\mathcal{G}) = \sum_{\beta \in \lspec(\mathcal{G})} \left|\beta - \frac{tr(D(\mathcal{G}))}{|V(\mathcal{G})|}\right|
\]
and
\[
LE^+(\mathcal{G}) = \sum_{\gamma \in \qspec(\mathcal{G})} \left| \gamma - \frac{tr(D(\mathcal{G}))}{|V(\mathcal{G})|}\right|,
\]
where $\spec(\mathcal{G}), \lspec(\mathcal{G})$ and $\qspec(\mathcal{G})$ are the spectrum, Laplacian spectrum and signless Laplacian spectrum of $\mathcal{G}$; $tr(D(\mathcal{G}))$ is the trace of $D(\mathcal{G})$ and  $V(\mathcal{G})$ is the set of vertices of $\mathcal{G}$.

In 1978, Gutman \cite{Gutman78} has introduced  the notion of $E({\mathcal{G}})$, which has been studied extensively by many mathematicians over the years (see \cite{Gutman-Furtula-2017} and the references therein).
In 2006, Gutman  and Zhou \cite{zhou} have introduced  the notion of $LE({\mathcal{G}})$; and in 2008, Gutman et al. \cite{GAVBR}   have introduced  the notion of $LE^{+}({\mathcal{G}})$.

 Let $N_{\mathcal{G}}(v_i)$ be the neighborhood set of vertex $v_i$ in $\mathcal{G}$. Also let $d_{\mathcal{G}}(v_i)$ be the degree of the vertex $v_i$ in $\mathcal{G}$, that is, $d_{\mathcal{G}}(v_i)=|N_{\mathcal{G}}(v_i)|$. Let $m_{\mathcal{G}}(v_i)$ be the average degree of the adjacent vertices of vertex $v_i$ in $\mathcal{G}$. If $v_i$ is an isolated vertex in $\mathcal{G}$, then we assume that $m_{\mathcal{G}}(v_i)=0$. Hence $d_{\mathcal{G}}(v_i)\,m_{\mathcal{G}}(v_i)=\sum\limits_{v_j:\,v_iv_j\in E(G)}\,d_{\mathcal{G}}(v_j)$. The common neighborhood of two vertices $v_i$ and $v_j$ is the set $N_{\mathcal{G}}(v_i)\cap N_{\mathcal{G}}(v_j)=N(v_i,\,v_j)$, (say), containing all the vertices other than $v_i$ and $v_j$ that are adjacent to both  $v_i$ and $v_j$. The common neighborhood matrix $\CN(\mathcal{G})$ of $\mathcal{G}$ is given by
\[
\CN(\mathcal{G})_{i, j} = \begin{cases}
|N(v_i, v_j)|  &\text{ if } i \ne j,\\[1mm]
0 &\text{ if } i = j,
\end{cases}
\]
where $\CN(\mathcal{G})_{i, j}$ is the $(i, j)$-th entry of $\CN(\mathcal{G})$. The common neighborhood spectrum (also known as CN-spectrum) of $\mathcal{G}$, denoted by $\cnspec(\mathcal{G})$,  is the set of all eigenvalues of  $\CN(\mathcal{G})$ with multiplicity.
We write $\cnspec(\mathcal{G}) = \{\mu_1^{a_1}, \mu_2^{a_2}, \dots,  \mu_k^{a_k}\}$, where the exponents $a_i$ of $\mu_i$ are the multiplicities of $\mu_i$ for $i = 1, 2, \dots ,k$; and $\mu_1, \mu_2, \dots, \mu_k$ are the distinct eigenvalues  of $\CN(\mathcal{G})$. The common neighborhood energy (also known as CN-energy) of $\mathcal{G}$, denoted by $E_{CN}(\mathcal{G})$, is defined as
\[
E_{CN}(\mathcal{G}) = \sum_{\mu \in \cnspec(\mathcal{G})} |\mu|.
\]
The concepts of CN-spectrum and  CN-energy are relatively new and not much explored. A graph $\mathcal{G}$ is called CN-hyperenergetic if $E_{CN}(\mathcal{G}) > E_{CN}(K_{|V(\mathcal{G})|})$, where $K_{|V(\mathcal{G})|}$ is the complete graph on $|V(\mathcal{G})|$-vertices. The concepts of CN-energy and CN-hyperenergetic graphs were introduced by Alwardi et al. \cite{ASG} in 2011. It is worth mentioning that the concept of hyperenergetic graph was introduced by Walikar et al. \cite{Walikar} and Gutman  \cite{Gutman1999}, independently in 1999.
Recall that $\mathcal{G}$ is called hyperenergetic if $E(\mathcal{G}) > E(K_{|V(\mathcal{G})|})$. L-hyperenergetic and Q-hyperenergetic graphs were defined in a similar way and these were introduced in \cite{FSN-2020}. Very recently, certain mathematical properties on $\CN(\mathcal{G})$ and $E_{CN}(\mathcal{G})$ of $\mathcal{G}$ are studied, see \cite{NFD} and the references cited therein.

In this paper we introduce the concepts of common neighborhood Laplacian spectrum, common neighborhood signless Laplacian spectrum and their corresponding energies of a graph
$\mathcal{G}$. In Section 2, we  establish relations between these energies and the first Zagreb index of a graph. In Section 3, we  introduce the  concepts of CNL-hyperenergetic and CNSL-hyperenergetic graphs and show that $\mathcal{G}$   is neither CNL-hyperenergetic nor CNSL-hyperenergetic if $\mathcal{G}$ is a complete bipartite graph.
In Section 4, we obtain certain relations between  various energies  of a graph.  Finally, in Section 5, we conclude this paper with several bounds for common neighborhood Laplacian and signless Laplacian energies of a graph.

\section{Definition and connection with Zagreb index}

First we observe that the $(i, j)$-th entry of $D(\mathcal{G})$ is given by
\[
D(\mathcal{G})_{i, j} = \begin{cases}
\underset{k = 1}{\overset{|V(\mathcal{G})|}{\sum}}A(\mathcal{G})_{i, k} &   \text{ if }  i = j  \text{ and } i = 1, 2, \dots, |V(\mathcal{G})|,\\[3mm]
0 & \text{ if } i \ne j,
\end{cases}
\]
where $A(\mathcal{G})_{i, k}$ is the $(i, k)$-th entry of $A(\mathcal{G})$. Thus $D(\mathcal{G})$ is a diagonal matrix whose diagonal entries are the corresponding row sums of the adjacency matrix of $\mathcal{G}$. Similarly, we define common neighborhood row sum matrix (abbreviated as CNRS-matrix) of $\mathcal{G}$ as given below:
\[
\CNRS(\mathcal{G})_{i, j} =  \begin{cases}
\underset{k = 1}{\overset{|V(\mathcal{G})|}{\sum}}\CN(\mathcal{G})_{i, k} &   \text{ if }  i = j  \text{ and } i = 1, 2, \dots, |V(\mathcal{G})|,\\[3mm]
0 & \text{ if } i \ne j,
\end{cases}
\]
where $\CNRS(\mathcal{G})$ is the CNRS-matrix of $\mathcal{G}$ and $\CNRS(\mathcal{G})_{i, j}$ is the $(i, j)$-th entry of $\CNRS(\mathcal{G})$.  The
common neighborhood Laplacian matrix and the common neighborhood signless Laplacian matrix (abbreviated as CNL-matrix and CNSL-matrix) of  $\mathcal{G}$, denoted by $\CNL(\mathcal{G})$ and $\CNSL(\mathcal{G})$, respectively, are defined as
\[
\CNL(\mathcal{G}):= \CNRS(\mathcal{G}) - \CN(\mathcal{G}) \text{ and } \CNSL(\mathcal{G}):= \CNRS(\mathcal{G}) + \CN(\mathcal{G}).
\]
Note that the matrices $\CNL(\mathcal{G})$ and $\CNSL(\mathcal{G})$ are symmetric and positive semidefinite.
The set of eigenvalues of  $\CNL(\mathcal{G})$ and $\CNSL(\mathcal{G})$  with multiplicities are called common neighborhood Laplacian spectrum and common neighborhood signless Laplacian spectrum (abbreviated as CN-Laplacian spectrum and CN-signless Laplacian spectrum) of $\mathcal{G}$, respectively. We write $\cnlspec(\mathcal{G})$ and $\cnqspec(\mathcal{G})$ to denote CN-Laplacian spectrum and CN-signless Laplacian spectrum of $\mathcal{G}$, respectively. By writing  $\cnlspec(\mathcal{G}) = \{\nu_1^{b_1}, \,\nu_2^{b_2}, \dots,\,\nu_{\ell}^{b_{\ell}}\}$ and $\cnqspec(\mathcal{G})= \{\sigma_1^{c_1},\,\sigma_2^{c_2}, \dots,\,\sigma_m^{c_m}\}$ we mean that $\nu_1,\,\nu_2, \dots, \,\nu_{\ell}$ are the distinct eigenvalues of $\CNL(\mathcal{G})$ with corresponding multiplicities  $b_1,\,b_2, \dots,\,b_{\ell}$ and $\sigma_1,\,\sigma_2, \dots,\,\sigma_m$ are the distinct eigenvalues of $\CNSL(\mathcal{G})$ with corresponding multiplicities  $c_1,\,c_2, \dots,\,c_m$.

Corresponding to CN-Laplacian spectrum and CN-signless Laplacian spectrum of $\mathcal{G}$ we define common neighborhood Laplacian energy and common neighborhood signless Laplacian energy (abbreviated as CNL-energy and CNSL-energy) of $\mathcal{G}$. The CNL-energy and CNSL-energy of $\mathcal{G}$, denoted by $LE_{CN}(\mathcal{G})$ and $LE^+_{CN}\mathcal{G}$, are as defined below:
\begin{equation}\label{LEcn}
LE_{CN}(\mathcal{G}) := \sum_{\nu \in \cnlspec(\mathcal{G})} \left|\nu - \frac{tr(\CNRS(\mathcal{G}))}{|V(\mathcal{G})|}\right|
\end{equation}
and
\begin{equation}\label{LE+cn}
LE^+_{CN}(\mathcal{G}) := \sum_{\sigma \in \cnqspec(\mathcal{G})} \left| \sigma - \frac{tr(\CNRS(\mathcal{G}))}{|V(\mathcal{G})|}\right|.
\end{equation}
This appends two new  entries in the list of energies prepared by Gutman and Furtula \cite{GUTMAN-2019*}. Following the concepts of various hyperenergetic graphs \cite{Walikar, Gutman1999, ASG, FSN-2020}, we introduce the concepts of CNL-hyperenergetic and CNSL-hyperenergetic graphs. A graph $\mathcal{G}$ is called CNL-hyperenergetic and CNSL-hyperenergetic if
\[
LE_{CN}(\mathcal{G}) > LE_{CN}(K_{|V(\mathcal{G})|}) \text{ and } LE^+_{CN}(\mathcal{G}) > LE^+_{CN}(K_{|V(\mathcal{G})|}),
\]
respectively.

The following lemma is useful in computing CN-Laplacian spectrum and   CN-signless Laplacian spectrum of a graph having disconnected components.
\begin{lem}\label{union-lem-1}
If $\mathcal{G} = \mathcal{G}_1 \sqcup \mathcal{G}_2 \sqcup \cdots \sqcup \mathcal{G}_k$ (that is, $\mathcal{G}_1, \mathcal{G}_2, \dots, \mathcal{G}_k$ are the disconnected components of $\mathcal{G}$), then
$$\cnlspec(\mathcal{G}) = \cnlspec(\mathcal{G}_1) \cup \cnlspec(\mathcal{G}_2) \cup \cdots \cup \cnlspec(\mathcal{G}_k)$$
and
$$\cnqspec(\mathcal{G}) = \cnqspec(\mathcal{G}_1) \cup \cnqspec(\mathcal{G}_2) \cup \cdots \cup \cnqspec(\mathcal{G}_k)$$
 counting multiplicities.
\end{lem}

Recall that first Zagreb index $M_1(\mathcal{G})$ is defined as
$$M_1(\mathcal{G})=\sum\limits^{|V(\mathcal{G})|}_{i=1}\,d_{\mathcal{G}}(v_i)^2=\sum\limits_{v_iv_j\in E(G)}\,\Big(d_{\mathcal{G}}(v_i)+d_{\mathcal{G}}(v_j)\Big).$$
Mathematical properties on the first Zagreb index was reported in \cite{BD1,DA1,DA2,DA3,DA4,DA5,DA6}. By \cite{DA4}, we have
\begin{equation}
	M_1(\mathcal{G})=\sum\limits^{|V(\mathcal{G})|}_{i=1}\,d_{\mathcal{G}}(v_i)^2=\sum\limits^{|V(\mathcal{G})|}_{i=1}\,d_{\mathcal{G}}(v_i)\,m_{\mathcal{G}}(v_i).\label{1kandas1}
\end{equation}
The following lemmas are useful in deriving relations between CN-Laplacian energy, CN-signless Laplacian energy and first Zagreb index.

\begin{lem}{\rm \cite{DA1}} \label{1kin1} Let $\mathcal{G}$ be a graph of order $|V(\mathcal{G})|$. Then for each $v_i\in V(\mathcal{G})$,
	$$\sum\limits^{|V(\mathcal{G})|}_{k=1,\,k\neq i}\,|N_{\mathcal{G}}(v_i)\cap N_{\mathcal{G}}(v_k)|=\sum\limits_{v_j:\,v_iv_j\in E(\mathcal{G})}\,\Big(d_{\mathcal{G}}(v_j)-1\Big)=d_{\mathcal{G}}(v_i)\,m_{\mathcal{G}}(v_i)-d_{\mathcal{G}}(v_i),$$
	where $d_{\mathcal{G}}(v_i)$ is the degree of the vertex $v_i$ in $\mathcal{G}$.
\end{lem}

\begin{lem}\label{1kd1} Let $\mathcal{G}$ be a graph with $|e(\mathcal{G})|$ edges and the first Zagreb index $M_1(\mathcal{G})$. Then $tr(\CNRS(\mathcal{G}))=M_1(\mathcal{G})-2\,|e(\mathcal{G})|$.
\end{lem}

\begin{proof} From the definition with Lemma \ref{1kin1} and (\ref{1kandas1}), we obtain
	\begin{align*}
		tr(\CNRS(\mathcal{G}))&=\sum\limits^{|V(\mathcal{G})|}_{i=1}\,\CNRS(\mathcal{G})_{i, i}\\[2mm]
		&=\sum\limits^{|V(\mathcal{G})|}_{i=1}\,\sum\limits^{|V(\mathcal{G})|}_{j=1,\,j\neq i}\,|N_{\mathcal{G}}(v_i)\cap N_{\mathcal{G}}(v_j)|\\[2mm]
		&=\sum\limits^{|V(\mathcal{G})|}_{i=1}\,\Big[d_{\mathcal{G}}(v_i)\,m_{\mathcal{G}}(v_i)-d_{\mathcal{G}}(v_i)\Big]\\[2mm]
		&=M_1(\mathcal{G})-2\,|e(\mathcal{G})|.
	\end{align*}
	This completes the result.
\end{proof}
We conclude this section with the following relations between $LE_{CN}(\mathcal{G})$, $	LE^+_{CN}(\mathcal{G})$ and   $M_1(\mathcal{G})$ which can be obtained from \eqref{LEcn}, \eqref{LE+cn} and  Lemma \ref{1kd1}.
\begin{thm}\label{relation_ZI}
Let $\mathcal{G}$ be a graph with $|e(\mathcal{G})|$ edges and the first Zagreb index $M_1(\mathcal{G})$. Then
\[
LE_{CN}(\mathcal{G}) = \sum_{\nu \in \cnlspec(\mathcal{G})} \left|\nu - \frac{M_1(\mathcal{G})-2\,|e(\mathcal{G})|}{|V(\mathcal{G})|}\right|
\]
and
\[
	LE^+_{CN}(\mathcal{G}) = \sum_{\sigma \in \cnqspec(\mathcal{G})} \left| \sigma - \frac{M_1(\mathcal{G})-2\,|e(\mathcal{G})|}{|V(\mathcal{G})|}\right|.
\]
\end{thm}

\section{CN-(signless) Laplacian spectrum and CN-(signless) Laplacian energy}
In this section we compute CN-(signless) Laplacian spectrum and CN-(signless) Laplacian energy of some classes of graphs and discuss their properties.
\begin{exam}\label{CN-LE-Kn}
For $n = 1$, it is clear that $\cnlspec(K_1) = \{0^1\}$, $\cnqspec(K_1) = \{0^1\}$  and so $LE_{CN} = 0$, $LE^+_{CN} = 0$. Therefore, we consider $n \geq 2$.
We have
$$\CN(K_n) = (n - 2)A(K_n)~\mbox{ and }~\CNRS(K_n) = \diag[(n-1)(n - 2), \dots,\,(n-1)(n-2)],$$
so $\CNL(K_n) = (n - 2)L(K_n)$ and $\CNSL(K_n) = (n - 2)Q(K_n)$. Also, $\lspec(K_n) = \{0^1, n^{n-1}\}$ and $\qspec(K_n) =\Big\{(2(n-1))^1, (n-2)^{n-1}\Big\}$. Therefore $\cnlspec(K_n) = \Big\{0^1, (n(n-2))^{n - 1}\Big\}$ and $\cnqspec(K_n) = \Big\{(2(n-1)(n-2))^1, ((n-2)^2)^{n - 1}\Big\}$.
We have $tr(\CNRS(K_n)) = n(n-1)(n - 2)$ and so $\frac{tr(\CNRS(K_n))}{|V(K_n)|} = (n-1)(n - 2)$. Therefore
\begin{align*}
&\left|0 -  \frac{tr(\CNRS(K_n))}{|V(K_n)|}\right|=(n-1)(n-2),\\[3mm]
&\left|n(n-2) -  \frac{tr(\CNRS(K_n))}{|V(K_n)|}\right|=n(n-2)-(n-1)(n - 2)=(n-2),
\end{align*}
and
\begin{align*}
&\left|2(n-1)(n-2) -  \frac{tr(\CNRS(K_n))}{|V(K_n)|}\right| = (n-1)(n-2),\\[3mm]
&\left|(n-2)^2 -  \frac{tr(\CNRS(K_n))}{|V(K_n)|}\right|=|(n-2)^2 - (n-1)(n - 2)|=|-(n-2)|= n-2.
\end{align*}
Hence, by \eqref{LEcn} and \eqref{LE+cn}, we obtain
\begin{align*}
&LE_{CN}(K_n) = (n-1)(n-2) + (n-1)(n-2) = 2(n-1)(n-2)\\
\mbox{and}&\\
&LE^+_{CN}(K_n) = (n-1)(n-2) + (n-1)(n-2) = 2(n-1)(n-2).
\end{align*}

By \cite[Proposition 2.4]{ASG}, it follows that if $\mathcal{G}_1$ and $\mathcal{G}_2$ are two disconnected component of $\mathcal{G}$, then $E_{CN}(\mathcal{G}) = E_{CN}(\mathcal{G}_1) + E_{CN}(\mathcal{G}_2)$. However, $LE_{CN}(\mathcal{G}) \ne LE_{CN}(\mathcal{G}_1) + LE_{CN}(\mathcal{G}_2)$ and $LE^+_{CN}(\mathcal{G}) \ne LE^+_{CN}(\mathcal{G}_1) + LE^+_{CN}(\mathcal{G}_2)$, if $\mathcal{G} = \mathcal{G}_1 \sqcup \mathcal{G}_2$. For example, if $\mathcal{G} = K_4 \sqcup K_6$ then, by Lemma \ref{union-lem-1} with the above result, it follows that
$$\cnlspec(\mathcal{G}) = \{0^2, 8^3, 24^5\}~\mbox{ and }~\cnqspec(\mathcal{G}) = \{4^3, (12)^1, (16)^5, (40)^1\}.$$
We have
\[
\frac{tr(\CNRS(\mathcal{G}))}{|V(\mathcal{G})|} = \frac{tr(\CNRS(K_4) + tr(\CNRS(K_6)}{10} = \frac{24 + 120}{10} = \frac{144}{10}.
\]
Therefore, by \eqref{LEcn} and \eqref{LE+cn}, we obtain
\begin{align*}
LE_{CN}(\mathcal{G}) &= 2\times\left|0 - \frac{144}{10}\right| + 3\times\left|8 - \frac{144}{10}\right| + 5\times\left|24 - \frac{144}{10}\right|\\
&= 2\times\frac{144}{10} + 3\times\frac{64}{10} + 5\times\frac{96}{10}= 96,
\end{align*}
but $LE_{CN}(K_4) + LE_{CN}(K_6) = 12 + 40 = 52$, and
\begin{align*}
LE^+_{CN}(\mathcal{G}) &= 3\times\left|4 - \frac{144}{10}\right| + 1\times\left|12 - \frac{144}{10}\right| + 5\times\left|16 - \frac{144}{10}\right| + 1\times\left|40 - \frac{144}{10}\right|\\
&= 2\times\frac{104}{10} + \frac{24}{10} + \frac{80}{10} + \frac{256}{10}= \frac{672}{10},
\end{align*}
but $LE^+_{CN}(K_4) + LE^+_{CN}(K_6) = 12 + 40 = 52$.
\end{exam}

\begin{exam}\label{CN-LE-Km,n} We now compute CN-(signless) Laplacian spectrum and CN-(signless) Laplacian energy of the complete bipartite graph $K_{m, n}$ on $(m + n)$-vertices. For this, let $V(K_{m, n}) \, = \,  \{v_1,\,v_2, \dots,\,v_m,\,v_{m+1},\,v_{m+2}, \dots,\,v_{m+n}\}$ \, and $\{v_1, v_2, \dots, v_m\}$, $\{v_{m+1}, v_{m+2}, \dots,$\\ $v_{m+n}\}$ be two partitions of $V(K_{m, n})$ such that every vertex in one set is adjacent to every vertex in the other set. We have
\[
\CN(K_{m,n}) =  \begin{bmatrix}
	n\,A(K_m) &  0 \\[2mm]
	0 & m\,A(K_n)
	\end{bmatrix}
\]
and
\[	
\CNRS(K_{m,n}) = \diag[\underset{m\text{-times}}{\underbrace{(m-1)n, \dots,(m-1)n}}, \underset{n\text{-times}}{\underbrace{(n-1)m, \dots, (n-1)m}}].
\]
Thus we have
$$\CNL(K_{m,n}) =  \begin{bmatrix}
	n\,L(K_m) &  0 \\[2mm]
	0 & m\,L(K_n)
	\end{bmatrix}
\mbox{ and }
\CNSL(K_{m,n}) =  \begin{bmatrix}
	n\,Q(K_m) &  0 \\[2mm]
	0 & m\,Q(K_n)
	\end{bmatrix}.$$

\vspace*{5mm}

\noindent
Since $\lspec(K_m) = \{0^1, m^{m-1}\}$ and $\lspec(K_n) = \{0^1, n^{n-1}\}$, therefore
$$\cnlspec(K_{m,n})= \Big\{(n\times 0)^1, (n\times m)^{m-1}, (m\times 0)^1, (m\times n)^{n-1}\Big\}= \Big\{0^2, (mn)^{m+n-2}\Big\}.$$
We have $\frac{tr(\CNRS(K_{m,n}))}{|V(K_{m,n})|} = \displaystyle{\frac{mn(m+n-2)}{m+n}}$ and so
\[
\left|0 -  \frac{tr(\CNRS(K_{m,n}))}{|V(K_{m,n})|}\right| = \frac{mn(m+n-2)}{m+n} \text{ and }
 \left|mn -  \frac{tr(\CNRS(K_{m,n}))}{|V(K_{m,n})|}\right| =  \frac{2mn}{m+n}.
\]

\vspace*{5mm}

\noindent
Hence, by \eqref{LEcn}, we get
$$LE_{CN}(K_{m, n}) = \frac{2mn(m+n-2)}{m+n} + \frac{2mn(m+n-2)}{m+n}  = \frac{4mn(m+n-2)}{m+n}.$$
In particular, for $m = n$, we obtain $\cnlspec(K_{m, n}) =\Big\{0^2,\,(m^2)^{2m-2}\Big\}$  and $LE_{CN}(K_{m, n}) = 4m(m-1)$.

\vspace*{3mm}

Again since $\qspec(K_m) = \Big\{(2(m-1))^1,\,(m-2)^{m-1}\Big\} \text{ and } \qspec(K_n) = \Big\{(2(n-1))^1,\,(n-2)^{n-1}\Big\}$, therefore,
$$\cnqspec(K_{m,n})= \Big\{(2n(m-1))^1,\,(n(m-2))^{m-1},\,(2m(n-1))^1,\,(m(n-2))^{n-1}\Big\}.$$
Note that if $m = n =1$, then $K_{1, 1} = K_2$. Hence, $\cnqspec(K_{m, n}) = \cnqspec(K_{2})$ $= \{0^2\}$ and $LE^+_{CN}(K_{1, 1}) = LE^+_{CN}(K_{2}) = 0$.
We now assume that $n\geq 2$ or $m\geq 2$. Now,
\begin{align*}
\left|2n(m-1) -  \frac{tr(\CNRS(K_{m,n}))}{|V(K_{m,n})|}\right| &= \left|2n(m-1) - \frac{mn(m+n-2)}{m+n}\right| \\[3mm]
&=\begin{cases}
\displaystyle{\frac{n(n-1)}{n+1}} & \text{ if } m = 1,\\[4mm]
\displaystyle{\frac{n(m+n)(m-2) + 2mn}{m+n}}  & \text{ if } m \geq 2,
\end{cases}
\end{align*}
\[
 \left|n(m-2) -  \frac{tr(\CNRS(K_{m,n}))}{|V(K_{m,n})|}\right| =  \left|n(m-2) - \frac{mn(m+n-2)}{m+n}\right| = \frac{2n^2}{m+n},
\]
\begin{align*}
\left|2m(n-1) -  \frac{tr(\CNRS(K_{m,n}))}{|V(K_{m,n})|}\right| &= \left|2m(n-1) - \frac{mn(m+n-2)}{m+n}\right| \\
&=\begin{cases}
\displaystyle{\frac{m(m-1)}{m+1}} & \text{ if } n = 1,\\[4mm]
\displaystyle{\frac{m(m+n)(n-2) + 2mn}{m+n}}  & \text{ if } n \geq 2,
\end{cases}
\end{align*}
and
\[
 \left|m(n-2) -  \frac{tr(\CNRS(K_{m,n}))}{|V(K_{m,n})|}\right| =  \left|m(n-2) - \frac{mn(m+n-2)}{m+n}\right| = \frac{2m^2}{m+n}.
\]
For $m=1$ and $n \geq 2$, by \eqref{LE+cn}, we have
\begin{align*}
LE^+_{CN}(K_{m, n}) &= \frac{n(n-1)}{n+1}
+ \frac{(n+1)(n-2) + 2n}{n+1} + \frac{2(n-1)}{n+1}= \frac{2(n-1)(n+2)}{n+1}.
\end{align*}
For $m \geq 2$ and $n = 1$, by \eqref{LE+cn}, we have
\begin{align*}
LE^+_{CN}(K_{m, n}) &= \frac{(m+1)(m-2) + 2m}{m+1} + \frac{2(m-1)}{m+1} + \frac{m(m-1)}{m+1}= \frac{2(m-1)(m+2)}{m+1}.
\end{align*}
For $m, n \geq 2$, by \eqref{LE+cn}, we have
\begin{align*}
&LE^+_{CN}(K_{m, n})\\[2mm]
=&\frac{n(m+n)(m-2) + 2mn}{m+n} + \frac{2n^2(m-1)}{m+n}+ \frac{m(m+n)(n-2) + 2mn}{m+n} + \frac{2m^2(n-1)}{m+n}\\[2mm]
=&\frac{4(m^2(n-1) + n^2(m-1))}{m+n}.
\end{align*}
Hence
\[
LE^+_{CN}(K_{m, n}) = \begin{cases}
0 & \text{ if } m = 1 \text{ and } n = 1,\\[2mm]
\displaystyle{\frac{2(n-1)(n+2)}{n+1}} & \text{ if } m = 1 \text{ and } n \geq 2, \\[4mm]
\displaystyle{\frac{2(m-1)(m+2)}{m+1}} & \text{ if } m \geq 2 \text{ and } n = 1, \\[4mm]
\displaystyle{\frac{4(m^2(n-1) + n^2(m-1))}{m+n}} & \text{ if } m, n \geq 2.
 \end{cases}
\]
In particular, for $m = n$, we obtain $\cnqspec(K_{m, n}) = \Big\{(2m(m-1))^2,\,(m(m-2))^{2(m-1)}\Big\}$
\text{ and }
$LE^+_{CN}(K_{m, n}) = 4m(m-1)$.
\end{exam}

\begin{prop}
Let $\overline{\mathcal{G}}$ be the complement of a graph $\mathcal{G}$ and $\mathcal{G}_1\vee \mathcal{G}_2$ be the join of two graphs $\mathcal{G}_1$ and $\mathcal{G}_2$.
\begin{enumerate}
\item If $\mathcal{G}=K_{n_1}\vee K_{n_2}\vee \cdots\vee K_{n_k}$, then
$$\cnlspec(\mathcal{G})=\{0^1,((n_1+n_2+\cdots+n_k)(n_1+n_2+\cdots+n_k-2))^{(n_1+n_2+\cdots+n_k-1)}\},$$
$$\!\!\!\!\!\!\!\!\!\!\!\!\cnqspec(\mathcal{G})=\{(2(n_1+n_2+\cdots+n_k-1)(n_1+n_2+\cdots+n_k-2))^1,((n_1+n_2+\cdots+n_k-2)^2)^{n-1}\}$$
and
$$LE_{CN}(\mathcal{G})=2(n_1+n_2+\cdots+n_k-1)(n_1+n_2+\cdots+n_k-2)=LE^+_{CN}(\mathcal{G}).$$
\item If $\mathcal{G}=\overline{K}_m\vee \overline{K}_n$, then
    \[
\cnlspec(\mathcal{G}) =\Big\{0^2, (mn)^{m + n-2}\Big\} \text{ and } LE_{CN}(\mathcal{G}) = \displaystyle{\frac{4mn(m+n-2)}{m+n}},
\]
\[
\cnqspec(\mathcal{G}) =\Big\{(2n(m-1))^1,(n(m-2))^{m-1},(2m(n-1))^1,(m(n-2))^{n-1}\Big\} \]
and
\[
LE^+_{CN}(\mathcal{G}) = \begin{cases}
\displaystyle{\frac{2(n-1)(n+2)}{n+1}} & \text{ if } m = 1 \text{ and } n \geq 2, \\[4mm]
\displaystyle{\frac{2(m-1)(m+2)}{m+1}} & \text{ if } m \geq 2 \text{ and } n = 1, \\[4mm]
\displaystyle{ \frac{4(m^2(n-1) + n^2(m-1))}{m+n}} & \text{ if } m, n \geq 2.
 \end{cases}
\]
\end{enumerate}
\end{prop}

\begin{proof}
The results follow from Examples {\ref{CN-LE-Kn}} and \ref{CN-LE-Km,n}  noting that $ K_{n_1}\vee K_{n_2}\vee \cdots\vee K_{n_k}=K_{n_1+n_2+\cdots+n_k}$ and $\overline{K_{m}}\vee \overline{K_{n}}=K_{m,n}$.	
\end{proof}

\begin{prop}\label{Kmn-CNL-nhyp}
The graph $K_{m, n}$ is not CNL-hyperenergetic.
\end{prop}

\begin{proof}
By Examples \ref{CN-LE-Kn} and \ref{CN-LE-Km,n}, we obtain
\begin{align*}
LE_{CN}(K_{m+n}) - LE_{CN}(K_{m, n}) &=  2(m+n-1)(m+n-2) -\frac{4mn(m+n-2)}{m+n}\\
&= \frac{[m(m-1) + n(n-1)]\,(m+n-2)}{m+n} \geq 0
\end{align*}
with equality if and only if $m=n=1$. Therefore,
\[
LE_{CN}(K_{m+n}) \geq LE_{CN}(K_{m, n})
\]
with equality if and only if $m=n=1$. Hence, the result follows.
\end{proof}

\begin{cor}
If $S_k$ denotes the star graph with one internal node and $k$ leaves, then
\[
\cnlspec(S_k) = \Big\{0^2,\, k^{k-1}\Big\} \text{ and } LE_{CN}(S_k) = \frac{4k(k-1)}{k+1}.
\]
Moreover, $S_k$ is not CNL-hyperenergetic.
\end{cor}

\begin{proof}
The result follows from Example \ref{CN-LE-Km,n} and Proposition \ref{Kmn-CNL-nhyp}, noting that $S_k = K_{1, k}$.
\end{proof}

\begin{prop}\label{Kmn-CNSL-nhyp}
The graph $K_{m, n}$ is not CNSL-hyperenergetic.
\end{prop}

\begin{proof}
If $m = n =1$ then we have $LE^+_{CN}(K_{1, 1}) = LE^+_{CN}(K_{2}) = 0$. Therefore, we consider the case when $m, n$ are not equal to $1$ simultaneously.
By Example \ref{CN-LE-Kn}, we have
\begin{equation}\label{LE+cn-Km+n}
LE^+_{CN}(K_{m+n}) = 2(m+n-1)(m+n-2).
\end{equation}
For $m=1$ and $n \geq 2$, by \eqref{LE+cn-Km+n} and Example \ref{CN-LE-Km,n}, we have
\begin{align*}
LE^+_{CN}(K_{m+n}) - LE^+_{CN}(K_{m, n}) &=   2n(n-1) - \frac{2(n-1)(n+2)}{n+1}\\[2mm]
&= \frac{2(n-1)(n^2-2)}{n+1} > 0.
\end{align*}
For $n=1$ and $m \geq 2$, by \eqref{LE+cn-Km+n} and Example \ref{CN-LE-Km,n}, we have
\begin{align*}
LE^+_{CN}(K_{m+n}) - LE^+_{CN}(K_{m, n}) &=   2m(m-1) - \frac{2(m-1)(m+2)}{m+1}\\[2mm]
&= \frac{2(m-1)(m^2-2)}{m+1} > 0.
\end{align*}
For $m, n \geq 2$, by \eqref{LE+cn-Km+n} and Example \ref{CN-LE-Km,n}, we have
\begin{align*}
LE^+_{CN}(K_{m+n}) &- LE^+_{CN}(K_{m, n}) \\
&=   2(m+n-1)(m+n-2) - \frac{4(m^2(n-1) + n^2(m-1))}{m+n}\\[2mm]
&= 2(m^3 + n^3 + m + n) + (m -n)^2 + 2mn(m + n - 2) > 0,
\end{align*}
since $m + n - 2 > 0$.   Therefore,
\[
LE^+_{CN}(K_{m+n}) \geq LE^+_{CN}(K_{m, n}).
\]
Hence, the result follows.
\end{proof}

\begin{cor}
If $S_k$ denotes the star graph with one internal node and $k$ leaves then
\[
\cnqspec(S_k) = \Big\{0^1, (2(k-1))^1, ((k-2))^{k-1}\Big\} \text{ and } LE^+_{CN}(S_k) = \frac{2(k-1)(k+2)}{k+1}.
\]
Moreover, $S_k$ is not CNSL-hyperenergetic.
\end{cor}

\begin{proof}
The result follows from Example \ref{CN-LE-Km,n} and Proposition \ref{Kmn-CNSL-nhyp}, noting that $S_k = K_{1, k}$.
\end{proof}

\section{Relation between various energies}
 In this section we derive some relations between $E_{CN}$, $LE_{CN}$, $LE^+_{CN}$, $E$, $LE$ and $LE^+$ of a graph $\mathcal{G}$. Let $e(\mathcal{G})$ be the set of edges of a graph $\mathcal{G}$.
Let $\spec(M)$ and  $E(M)$ be the spectrum and  energy of a real square symmetric matrix $M$ of size $n$. Then $\spec(M) = \{\mu_1,\mu_2,\ldots,\mu_n\}$,   where $\mu_1, \mu_2, \ldots, \mu_n$ are eigenvalues (not necessarily distinct) of $M$, and
\[
E(M)=\sum\limits_{i=1}^{n}|\mu_i|.
\]
The following two lemmas are useful in obtaining a relation between $E_{CN}(\mathcal{G})$ and $E(\mathcal{G})$.

\begin{lem}{\rm \cite[Proposition 2.7]{ASG}}\label{AG-DG-relation}
	Let $\mathcal{G}$ be any graph. Then $\CN(\mathcal{G})=A(\mathcal{G})^2-D(\mathcal{G})$.
\end{lem}

\begin{lem}{\rm \cite[Lemma 2.10]{DMG}}\label{Energy-inequality}
Let $M_1$ and $M_2$ be two real square symmetric matrices of order $n$ and let $M=M_1+M_2$. Then
$$E(M)\leq E(M_1)+E(M_2).$$
\end{lem}

\begin{thm}\label{CN energy inequality}
Let $\mathcal{G}$ be any graph with $|e(\mathcal{G})|$ edges. Then $E_{CN}(\mathcal{G})\leq E(\mathcal{G})^2+2 |e(\mathcal{G})|$.
\end{thm}

\begin{proof}
By  Lemmas \ref{AG-DG-relation} and \ref{Energy-inequality}, we obtain
\begin{align}\label{E-ECN-eq-1}
    E_{CN}(\mathcal{G})=&E(CN(\mathcal{G}))\nonumber \\
    =&E(A(\mathcal{G})^2-D(\mathcal{G}))\nonumber\\
    =&E(A(\mathcal{G})^2+(-D(\mathcal{G})))
    \leq E(A(\mathcal{G})^2)+E(-D(\mathcal{G})).
\end{align}
Let $\spec(\mathcal{G})=\{\lambda_1,\lambda_2,\ldots,\lambda_n\}$, where $n = |V(\mathcal{G})|$. Then $\spec(A(\mathcal{G})^2)=\{\lambda_1^2,\lambda_2^2,\ldots,\lambda_n^2\}$. Therefore,
\begin{align*}
    E(A(\mathcal{G})^2)=&\sum\limits_{i=1}^n|\lambda_i^2|
                        \leq \left(\sum\limits_{i=1}^n|\lambda_i|\right)^2
                        =E(\mathcal{G})^2.
\end{align*}
Again, let $\spec(D(\mathcal{G}))=\{d_{\mathcal{G}}(v_1), d_{\mathcal{G}}(v_2), \ldots, d_{\mathcal{G}}(v_n)\}$. Then
$\spec(-D(\mathcal{G}))=\{-d_{\mathcal{G}}(v_1),$ $-d_{\mathcal{G}}(v_2),\ldots,-d_{\mathcal{G}}(v_n)\}$. Therefore
$$E(-D(\mathcal{G}))=\sum\limits_{i=1}^n|-d_{\mathcal{G}}(v_i)| = \sum\limits_{i=1}^n d_{\mathcal{G}}(v_i) = 2\times |e(\mathcal{G})| = E(D(\mathcal{G})).$$
Hence, the result follows from \eqref{E-ECN-eq-1}.
\end{proof}

\begin{cor} Let $\mathcal{G}$ be any graph of order $|V(\mathcal{G})|$ with $|e(\mathcal{G})|$ edges. Then $E_{CN}(\mathcal{G})\leq 2 |e(\mathcal{G})|\,(|V(\mathcal{G})|+1)$.
\end{cor}

\begin{proof} It is well-known that $E(\mathcal{G})\leq \sqrt{2\,|e(\mathcal{G})|\,|V(\mathcal{G})|}$.
Using the above result with Theorem \ref{CN energy inequality}, we obtain
\begin{align*}
E_{CN}(\mathcal{G})\leq E(\mathcal{G})^2+2 |e(\mathcal{G})|
  \leq 2\,|e(\mathcal{G})|\,|V(\mathcal{G})|+2 |e(\mathcal{G})|.
\end{align*}
Hence the result follows.
\end{proof}

\begin{thm}\label{CNL and CNSL energy inequality}
Let $\mathcal{G}$ be any graph with $|e(\mathcal{G})|$ edges and the first Zagreb index $M_1(\mathcal{G})$. Then $LE_{CN}(\mathcal{G})\leq E_{CN}(\mathcal{G})+2\,\Big(M_1(\mathcal{G})-2\,|e(\mathcal{G})|\Big)$ and $LE^+_{CN}(\mathcal{G})\leq E_{CN}(\mathcal{G})+2\,\Big(M_1(\mathcal{G})-2\,|e(\mathcal{G})|\Big)$.
\end{thm}

\begin{proof}
	By Theorem  \ref{relation_ZI}, we have
	\begin{align*}
		LE_{CN}(\mathcal{G}) =& \sum_{\nu \in \cnlspec(\mathcal{G})} \left|\nu - \frac{M_1(\mathcal{G})-2\,|e(\mathcal{G})|}{|V(\mathcal{G})|}\right| \\
		\leq& \sum_{\nu \in \cnlspec(\mathcal{G})}|\nu|+\sum_{\nu \in \cnlspec(\mathcal{G})}\left|\frac{M_1(\mathcal{G})-2\,|e(\mathcal{G})|}{|V(\mathcal{G})|}\right|\\
		=& E(\CNL(\mathcal{G}))+\frac{M_1(\mathcal{G})-2\,|e(\mathcal{G})|}{|V(\mathcal{G})|}\sum_{\nu \in \cnlspec(\mathcal{G})}1\\
		=& E(\CNRS(\mathcal{G})-\CN(\mathcal{G}))+ M_1(\mathcal{G})-2\,|e(\mathcal{G})|.
	\end{align*}
	Using Lemmas \ref{Energy-inequality} \& \ref{1kd1} and the fact that  $E(\CN(\mathcal{G})) = E(-\CN(\mathcal{G}))$ \& $ E(\CNRS(\mathcal{G})) = tr(\CNRS(\mathcal{G}))$, we obtain
	\begin{align*}
		LE_{CN}(\mathcal{G}) &\leq E(\CNRS(\mathcal{G}))+E(\CN(\mathcal{G}))+M_1(\mathcal{G})-2\,|e(\mathcal{G})|\\
		&=E_{CN}(\mathcal{G})+2\,\Big(M_1(\mathcal{G})-2\,|e(\mathcal{G})|\Big).
	\end{align*}
	Similarly, the  bound for $LE^+_{CN}(\mathcal{G})$ follows from Theorem  \ref{relation_ZI} and Lemma \ref{Energy-inequality}.
\end{proof}
As a consequence of  Theorems \ref{CN energy inequality} and \ref{CNL and CNSL energy inequality}, we get the following relations between $E(\mathcal{G})$, $LE_{CN}(\mathcal{G})$ and $LE^+_{CN}(\mathcal{G})$.
\begin{cor}\label{bound}
Let $\mathcal{G}$ be any graph with $|e(\mathcal{G})|$ edges and the first Zagreb index $M_1(\mathcal{G})$. Then $LE_{CN}(\mathcal{G})$ and $LE^+_{CN}(\mathcal{G})$ are bounded above by
$$E(\mathcal{G})^2+2\,\Big(M_1(\mathcal{G})-|e(\mathcal{G})|\Big).$$
\end{cor}

\begin{rem}\label{R1}
Theorem \ref{CN energy inequality} gives relation between $E_{CN}(\mathcal{G})$ and $E(\mathcal{G})$. Theorem \ref{CNL and CNSL energy inequality} and Corollary \ref{bound} give relations between $LE_{CN}(\mathcal{G})$, $LE^+_{CN}(\mathcal{G})$, $E_{CN}(\mathcal{G})$ and $E(\mathcal{G})$. However, using the facts that
\[
E(A(\mathcal{G})^2)=\sum\limits_{i=1}^n|\lambda_i^2|=2|e(\mathcal{G})|\leq E(\mathcal{G})^2
\]
and
\[
\sum\limits_{\nu\in\cnlspec(\mathcal{G})}|\nu|=tr(\CNRS(\mathcal{G}))=M_1(\mathcal{G})-2\,|e(\mathcal{G})|=\sum\limits_{\sigma\in\cnqspec(\mathcal{G})}|\sigma|,
\]
we get the following better upper bounds for $E_{CN}(\mathcal{G})$, $LE_{CN}(\mathcal{G})$ and $LE^+_{CN}(\mathcal{G})$:
\begin{equation}\label{rem-eq-1}
E_{CN}(\mathcal{G})\leq 4|e(\mathcal{G})|\leq 2E(\mathcal{G})^2,
\end{equation}
\begin{equation}\label{rem-eq-2}
LE_{CN}(\mathcal{G})\leq 2\,tr(\CNRS(\mathcal{G}))=2\,\Big(M_1(\mathcal{G})-2\,|e(\mathcal{G})|\Big)\geq LE^+_{CN}(\mathcal{G}).
\end{equation}
In Section 5, we shall obtain more bounds for $LE_{CN}(\mathcal{G})$ and $LE^+_{CN}(\mathcal{G})$.
\end{rem}
Recall that the derived graph of $\mathcal{G}$, denoted by $\mathcal{G}^{\dagger}$ is the graph with vertex set $V(\mathcal{G})$, in which two vertices are adjacent if and only if their distance in $\mathcal{G}$ is two.
\begin{thm}\label{Tri-Quad}
	If $\mathcal{G}$ is a triangle- and quadrangle-free graph, then $LE_{CN}(\mathcal{G})=LE(\mathcal{\mathcal{G}}^{\dagger})$ and $LE^+_{CN}(\mathcal{G})=LE^+(\mathcal{\mathcal{G}}^{\dagger})$, where $\mathcal{G}^{\dagger}$ is the derived graph of $\mathcal{G}$.
\end{thm}
\begin{proof}
	If $\mathcal{G}$ is a triangle- and quadrangle-free graph, then $CN(\mathcal{G})=A(\mathcal{G}^\dagger)$. Therefore, $\CNRS(\mathcal{G})=D(\mathcal{G}^\dagger)$ and so  $\CNL(G)=L(\mathcal{G}^\dagger)$ and $\CNSL(\mathcal{G})=Q(\mathcal{G}^\dagger)$. Hence, $\cnlspec(\mathcal{G}) = \lspec(\mathcal{G}^\dagger)$ and $\cnqspec(\mathcal{G}) = \qspec(\mathcal{G}^\dagger)$. Since $tr(\CNRS(\mathcal{G})) = tr(D(\mathcal{G}^\dagger))$ and $V(\mathcal{G}) = V(\mathcal{G}^\dagger)$, by \eqref{LEcn} and \eqref{LE+cn}, we have
	\[
	LE_{CN}(\mathcal{G}) = \sum_{\nu \in \lspec(\mathcal{G}^\dagger)} \left|\nu - \frac{tr(D(\mathcal{G}^\dagger))}{|V(\mathcal{G})|}\right| = LE(\mathcal{G}^{\dagger})
	\]
	and
	\[
	LE^+_{CN}(\mathcal{G})= \sum_{\sigma \in \qspec(\mathcal{G}^\dagger)} \left| \sigma - \frac{tr(D(\mathcal{G}^\dagger))}{|V(\mathcal{G})|}\right| =LE^+(\mathcal{G}^{\dagger}).
	\]
\end{proof}

\begin{cor}
	\begin{enumerate}	
		\item If $T$ is a tree then $LE_{CN}(T)=LE(T^{\dagger})\text{ and } LE^+_{CN}(T)=LE^+(T^{\dagger}).$
		\item If $P_n$ is the path on  $n$ vertices, then
		$LE_{CN}(P_n)=LE(P_{\lceil \frac{n}{2}\rceil})+LE(P_{\lfloor \frac{n}{2}\rfloor})$  and  $LE^+_{CN}(P_n)=LE^+(P_{\lceil \frac{n}{2}\rceil})+LE^+(P_{\lfloor \frac{n}{2}\rfloor}).$
		\item Let $C_n$ be the  cycle on  $n$ vertices.
		\begin{enumerate}
			\item If $n$ is odd and $n\geq 3$, then $LE_{CN}(C_n)=LE(C_n)$ and $LE^+_{CN}(C_n)=LE^+(C_n)$.
			\item If $n$ is even and $n> 4$, then $LE_{CN}(C_n)=2LE(C_{\frac{n}{2}})$ and $LE^+_{CN}(C_n)=2LE^+(C_{\frac{n}{2}})$. Also, $LE_{CN}(C_4) = 2LE(C_4) = LE^+_{CN}(C_4) = 2LE^+(C_4) = 8$.
		\end{enumerate}
	\end{enumerate}
\end{cor}
\begin{proof}
	(a) Follows from Theorem \ref{Tri-Quad} noting that   $T$ is  triangle- and quadrangle-free.

	\noindent (b) Follows from Theorem \ref{Tri-Quad} noting that  $P_n$ is  triangle- and quadrangle-free and
	$$P_n^\dagger\cong P_{\lceil \frac{n}{2}\rceil}\cup P_{\lfloor \frac{n}{2}\rfloor}.$$
	
	\noindent (c) Let $C_n$ be the  cycle on  $n$ vertices.
	
	(i) If $n=3$, then $C_3 = K_3$. Therefore,  $LE_{CN}(C_3) = LE^+_{CN}(C_3) = 4 = LE(C_3) = LE^+(C_3)$.

	If $n$ is odd and $n>3$, then $C_n$ is triangle- and quadrangle-free. Also,  $(C_n)^\dagger\cong C_n$. Hence, the result follows from Theorem \ref{Tri-Quad}.
	
	(ii) If $n$ is even and $n > 4$, then $C_n$ is triangle- and quadrangle-free. Also, $(C_n)^\dagger\cong C_{\frac{n}{2}}\cup C_{\frac{n}{2}}$. Therefore, by  Theorem \ref{Tri-Quad} we get
	\[
	LE_{CN}(C_n) = LE(C_{\frac{n}{2}} \cup C_{\frac{n}{2}}) =  2LE(C_{\frac{n}{2}})
	\]
	and
	\[
	LE^+_{CN}(C_n)= LE^+(C_{\frac{n}{2}} \cup C_{\frac{n}{2}}) = 2LE^+(C_{\frac{n}{2}}).
	\]
	
	If $n = 4$ then it is easy to see that $\CNRS(C_4) = D(C_4)$, which is a $4 \times 4$ diagonal matrix such that every element in the diagonal is equal to $2$, and  $\cnlspec(C_4) = \cnqspec(C_4) = \{0^2, 4^2\}$. Therefore, by \eqref{LEcn} and \eqref{LE+cn}, we have
	\[
	LE_{CN}(C_4) = LE^+_{CN}(C_4) = 8.
	\]
	Again,  $\lspec(C_4) = \qspec(C_4) = \{0, 2^2, 4\}$ and so $LE(C_4) = LE^+(C_4) = 4$. Thus, $LE_{CN}(C_4) = 2LE(C_4)$ and $LE^+_{CN}(C_4) = 2LE^+(C_4)$.
\end{proof}

\section{More bounds for $LE_{CN}(\mathcal{G})$ and $LE^+_{CN}(\mathcal{G})$}

 In this section we shall obtain several bounds for  $LE_{CN}(\mathcal{G})$ and $LE^+_{CN}(\mathcal{G})$.
 Since the matrices $\CNL(\mathcal{G})$ and $\CNSL(\mathcal{G})$ are positive semidefinite, the elements of $\cnlspec(\mathcal{G})$ and $\cnqspec(\mathcal{G})$ are non-negative. Thus we may write $\cnlspec(\mathcal{G})=\{\nu_1,\nu_2,\dots,$ $\nu_{|V(\mathcal{G})|}\}$ and $\cnqspec(\mathcal{G})=\{\sigma_1, \sigma_2, \dots, \sigma_{|V(\mathcal{G})|}\}$, where $\nu_1\geq \nu_2\geq\dots\geq\nu_{|V(\mathcal{G})|}$ and $\sigma_1\geq \sigma_2\geq\dots\geq\sigma_{|V(\mathcal{G})|}$. We have
 $$\sum\limits_{i=1}^{|V(\mathcal{G})|}\nu_i=\sum\limits_{\nu\in\cnlspec{(\mathcal{G})}}\nu=tr(\CNRS(\mathcal{G}))=\sum\limits_{\sigma\in\cnqspec{(\mathcal{G})}}\sigma=\sum\limits_{i=1}^{|V(\mathcal{G})|}\sigma_i.$$
 Also,\\
 $$\sum\limits_{\nu\in\cnlspec(\mathcal{G})}\left(\nu-\frac{tr(\CNRS(\mathcal{G}))}{|V(\mathcal{G})|}\right)=\sum\limits_{\sigma\in\cnqspec(\mathcal{G})}\left(\sigma-\frac{tr(\CNRS(\mathcal{G}))}{|V(\mathcal{G})|}\right)=0.$$

 \vspace*{3mm}

 \noindent
 Let $\alpha,\beta~(1\leq \alpha,\beta\leq |V(\mathcal{G})|)$ be the largest integers such that
 \begin{equation}\label{alpha-beta}
     \nu_\alpha\geq\frac{tr(\CNRS(\mathcal{G}))}{|V(\mathcal{G})|}=\frac{M_1(\mathcal{G})-2\,|e(\mathcal{G})|}{|V(\mathcal{G})|} \text{ and } \sigma_\beta\geq \frac{tr(\CNRS(\mathcal{G}))}{|V(\mathcal{G})|}= \frac{M_1(\mathcal{G})-2\,|e(\mathcal{G})|}{|V(\mathcal{G})|}.
 \end{equation}
Let  $S_\alpha(\mathcal{G})=\sum\limits_{i=1}^\alpha\nu_i$ and $S^+_\beta(\mathcal{G})=\sum\limits_{i=1}^\beta\sigma_i$. Then we have the  following useful lemmas.

\begin{lem}\label{Finite-sum equality}
For any graph $\mathcal{G}$, we have
\[LE_{CN}(\mathcal{G})=2S_\alpha(\mathcal{G})-\frac{2\alpha\,tr(\CNRS(\mathcal{G}))}{|V(\mathcal{G})|}=2S_\alpha(\mathcal{G})-\frac{2\alpha\,\Big(M_1(\mathcal{G})-2\,|e(\mathcal{G})|\Big)}{|V(\mathcal{G})|}\]
and
\[ LE^+_{CN}(\mathcal{G})=2S^+_\beta(\mathcal{G})-\frac{2\beta\,tr(\CNRS(\mathcal{G}))}{|V(\mathcal{G})|}=2S^+_\beta(\mathcal{G})-\frac{2\beta\,\Big(M_1(\mathcal{G})-2\,|e(\mathcal{G})|\Big)}{|V(\mathcal{G})|},\]
where $|e(\mathcal{G})|$ is the number of edges and $M_1(\mathcal{G})$ is the first Zagreb index in $\mathcal{G}$.
\end{lem}

\begin{lem}\label{maximum}
For any graph $\mathcal{G}$, we have
\[
LE_{CN}(\mathcal{G})=\max\limits_{1\leq i\leq |V(\mathcal{G})|}\left\lbrace\,2S_i(\mathcal{G})-\frac{2i\,\Big(M_1(\mathcal{G})-2\,|e(\mathcal{G})|\Big)}{|V(\mathcal{G})|}\right\rbrace
\]
and
\[
LE^+_{CN}(\mathcal{G})=\max\limits_{1\leq i\leq |V(\mathcal{G})|}\left\lbrace 2S^+_i(\mathcal{G})-\frac{2i\,\Big(M_1(\mathcal{G})-2\,|e(\mathcal{G})|\Big)}{|V(\mathcal{G})|}\right\rbrace,
\]
where $|e(\mathcal{G})|$ is the number of edges and $M_1(\mathcal{G})$ is the first Zagreb index in $\mathcal{G}$.
\end{lem}

\begin{proof}
Let $k~(1\leq k\leq |V(\mathcal{G})|)$ be any integer. For $k<\alpha$, by (\ref{alpha-beta}), we obtain
\[
S_\alpha(\mathcal{G})-S_k(\mathcal{G})=\sum\limits_{i=k+1}^\alpha\nu_i\geq \frac{(\alpha-k)\,\Big(M_1(\mathcal{G})-2\,|e(\mathcal{G})|\Big)}{|V(\mathcal{G})|}.
\]
For $k>\alpha $, we obtain
\[
S_k(\mathcal{G})-S_\alpha (\mathcal{G})=\sum\limits_{i=\alpha +1}^k\nu_i< \frac{(k-\alpha)\,\Big(M_1(\mathcal{G})-2\,|e(\mathcal{G})|\Big)}{|V(\mathcal{G})|},
\]
that is,
\[
S_\alpha(\mathcal{G})-S_k(\mathcal{G})> \frac{(\alpha-k)\,\Big(M_1(\mathcal{G})-2\,|e(\mathcal{G})|\Big)}{|V(\mathcal{G})|}.
\]
Moreover, $S_\alpha(\mathcal{G})=S_k(\mathcal{G})$ for $k=\alpha$. Thus for any value of $k~(1\leq k\leq |V(\mathcal{G})|)$, we obtain
\begin{align*}
    &S_\alpha(\mathcal{G})-S_k(\mathcal{G})\geq \frac{(\alpha-k)\,\Big(M_1(\mathcal{G})-2\,|e(\mathcal{G})|\Big)}{|V(\mathcal{G})|}
\end{align*}
and so
\[
2S_\alpha(\mathcal{G})-\frac{2\alpha\, tr(\CNRS(\mathcal{G}))}{|V(\mathcal{G})|}\geq 2S_k(\mathcal{G})-\frac{2k\,\Big(M_1(\mathcal{G})-2\,|e(\mathcal{G})|\Big)}{|V(\mathcal{G})|}.
\]
This gives
\[
2S_\alpha(\mathcal{G})-\frac{2\alpha\, tr(\CNRS(\mathcal{G}))}{|V(\mathcal{G})|}=\max\limits_{1\leq i\leq |V(\mathcal{G})|} \left\lbrace 2S_i(\mathcal{G})-\frac{2i\, \Big(M_1(\mathcal{G})-2\,|e(\mathcal{G})|\Big)}{|V(\mathcal{G})|}\right\rbrace.
\]
Similarly, it can be seen that
\[
2S^+_\beta(\mathcal{G})-\frac{2\beta\,tr(\CNRS(\mathcal{G}))}{|V(\mathcal{G})|}=\max\limits_{1\leq i\leq |V(\mathcal{G})|} \left\lbrace 2S^+_i(\mathcal{G})-\frac{2i\, \Big(M_1(\mathcal{G})-2\,|e(\mathcal{G})|\Big)}{|V(\mathcal{G})|}\right\rbrace.
\]
Hence, the result follows from  Lemma \ref{Finite-sum equality}.
\end{proof}

Let $(a) := (a_1, a_2, \cdots, a_n) \in \mathbb{R}^n$ and $(b):= (b_1, b_2, \cdots, b_n) \in \mathbb{R}^n$ be such that
$a_1 \geq a_2\geq\cdots\geq a_n$ and $b_1\geq b_2\geq\cdots\geq b_n$. Then $(a)$ is said to be majorize $(b)$ if
\[
\sum\limits_{i=1}^{k}a_i\geq\sum\limits_{i=1}^{k}b_i \text{ for } 1\leq k\leq n-1~~\text{ and }~~\sum\limits_{i=1}^{n}a_i=\sum\limits_{i=1}^{n}b_i.
\]
It is well-known that the spectrum of any symmetric, positive semidefinite matrix majorizes its main diagonal (see \cite{Schur}, \cite[p. 218]{MO79} as noted in \cite{GM-1994}). Since $\CNL(\mathcal{G})$ and $\CNSL(\mathcal{G})$ are symmetric and positive semidefinite  for any graph $\mathcal{G}$, we have the following lemma when the elements of $\cnlspec(\mathcal{G})$, $\cnqspec(\mathcal{G})$ and main diagonal elements of $\CNRS(\mathcal{G})$ are arranged in decreasing order.

\begin{lem}\label{majorize}
For any graph $\mathcal{G}$, $\cnlspec(\mathcal{G})$ and $\cnqspec(\mathcal{G})$  majorize main diagonal elements of $\CNRS(\mathcal{G})$ when the elements of $\cnlspec(\mathcal{G})$, $\cnqspec(\mathcal{G})$ and main diagonal elements of $\CNRS(\mathcal{G})$ are arranged in decreasing order.
\end{lem}
We write the main diagonal elements of $\CNRS(\mathcal{G})$ as  $\CNRS(\mathcal{G})_{i,i}$ for $1 \leq i \leq |V(\mathcal{G})|$, where $\CNRS(\mathcal{G})_{1,1} \geq \CNRS(\mathcal{G})_{2,2} \geq \cdots \geq \CNRS(\mathcal{G})_{|V(\mathcal{G})|,|V(\mathcal{G})|}$.
Now we give lower bounds for $LE_{CN}(\mathcal{G})$ and $LE^+_{CN}(\mathcal{G})$ analogous to the bound given by \cite[Theorem 3.1]{DM-2014} for $LE(\mathcal{G})$.

\begin{thm}\label{triangleleft-inequality}
Let $\mathcal{G}$ be a graph with $|e(\mathcal{G})|$ edges and the first Zagreb index $M_1(\mathcal{G})$. Then
\[
LE_{CN}(\mathcal{G})\geq 2\left(\Delta\,(\delta-1)-\frac{M_1(\mathcal{G})-2\,|e(\mathcal{G})|}{|V(\mathcal{G})|}\right)
\]
and
\[
LE^+_{CN}(\mathcal{G})\geq 2\left(\Delta\,(\delta-1)-\frac{M_1(\mathcal{G})-2\,|e(\mathcal{G})|}{|V(\mathcal{G})|}\right),
\]
where $\Delta$ and $\delta$ are the maximum degree and the minimum degree in $\mathcal{G}$, respectively.
\end{thm}

\begin{proof} Let $v_1$ be the maximum degree vertex in $\mathcal{G}$. Then $d_{\mathcal{G}}(v_1)=\Delta$ and $m_{\mathcal{G}}(v_1)\geq \delta$ as $\delta$ is the minimum degree in $\mathcal{G}$. As a consequence of Lemma \ref{majorize} with Lemma \ref{1kin1}, we obtain
\begin{align*}
    &\nu_1\geq \CNRS(\mathcal{G})_{1,1}=d_{\mathcal{G}}(v_1)\,m_{\mathcal{G}}(v_1)-d_{\mathcal{G}}(v_1)=\Delta\,\Big(m_{\mathcal{G}}(v_1)-1\Big)\geq \Delta\,(\delta-1) \\ \text{ and }&\\
& \sigma_1\geq \CNRS(\mathcal{G})_{1,1}=d_{\mathcal{G}}(v_1)\,m_{\mathcal{G}}(v_1)-d_{\mathcal{G}}(v_1)=\Delta\,\Big(m_{\mathcal{G}}(v_1)-1\Big)\geq \Delta\,(\delta-1),
\end{align*}
Using the above result with Lemma \ref{maximum}, we obtain
\begin{align*}
LE_{CN}(\mathcal{G}) \geq 2\,S_1(\mathcal{G})-\frac{2\,\Big(M_1(\mathcal{G})-2\,|e(\mathcal{G})|\Big)}{|V(\mathcal{G})|}&=2\, \nu_1-\frac{2\,\Big(M_1(\mathcal{G})-2\,|e(\mathcal{G})|\Big)}{|V(\mathcal{G})|}\\[3mm]
    &\geq 2\left(\Delta\,(\delta-1)-\frac{M_1(\mathcal{G})-2\,|e(\mathcal{G})|}{|V(\mathcal{G})|}\right).
\end{align*}
Similarly,
\begin{equation*}
LE^+_{CN}(\mathcal{G})\geq 2\, \sigma_1-\frac{2\,\Big(M_1(\mathcal{G})-2\,|e(\mathcal{G})|\Big)}{|V(\mathcal{G})|}\geq 2\left(\Delta\,(\delta-1)-\frac{M_1(\mathcal{G})-2\,|e(\mathcal{G})|}{|V(\mathcal{G})|}\right).
\end{equation*}
\end{proof}

\begin{thm}\label{diagonal-inequality} Let $\mathcal{G}$ be a graph with $|e(\mathcal{G})|$ edges and the first Zagreb index $M_1(\mathcal{G})$. Then
\[
LE_{CN}(\mathcal{G})\geq 2\left (\sum\limits_{i=1}^\alpha \CNRS(\mathcal{G})_{i,i}-\frac{\alpha\,\Big(M_1(\mathcal{G})-2\,|e(\mathcal{G})|\Big)}{|V(\mathcal{G})|} \right )
\]
and
\[
LE^+_{CN}(\mathcal{G})\geq 2\left (\sum\limits_{i=1}^\beta \CNRS(\mathcal{G})_{i,i}-\frac{\beta\,\Big(M_1(\mathcal{G})-2\,|e(\mathcal{G})|\Big)}{|V(\mathcal{G})|} \right ),
\]
where $\alpha$ and $\beta$ are as given in \eqref{alpha-beta} and $\CNRS(\mathcal{G})_{i,i}=d_{\mathcal{G}}(v_i)\,(m_{\mathcal{G}}(v_i)-1)$.
\end{thm}

\begin{proof} By Lemma \ref{1kin1}, we have $\CNRS(\mathcal{G})_{i,i}=d_{\mathcal{G}}(v_i)\,(m_{\mathcal{G}}(v_i)-1)$. By Lemma \ref{majorize}, we obtain
\[
\sum\limits_{i=1}^k\nu_i\geq \sum\limits_{i=1}^k \CNRS(\mathcal{G})_{i,i}\text{ and } \sum\limits_{i=1}^k\sigma_i\geq \sum\limits_{i=1}^k \CNRS(\mathcal{G})_{i,i}\, \text{ for } 1\leq k\leq |V(\mathcal{G})|.
\]
In particular, we have
\[
\sum\limits_{i=1}^\alpha\nu_i\geq \sum\limits_{i=1}^\alpha \CNRS(\mathcal{G})_{i,i}\text{ and } \sum\limits_{i=1}^\beta\sigma_i\geq \sum\limits_{i=1}^\beta \CNRS(\mathcal{G})_{i,i}.
\]
Therefore,
\[
S_\alpha(\mathcal{G})\geq \sum\limits_{i=1}^\alpha \CNRS(\mathcal{G})_{i,i} \text{ and } S^+_\beta(\mathcal{G})\geq \sum\limits_{i=1}^\beta \CNRS(\mathcal{G})_{i,i}.
\]
Hence, the result follows from Lemma \ref{Finite-sum equality}.
\end{proof}

Using Lemma \ref{Finite-sum equality}, we also have the following upper bounds for $LE_{CN}(\mathcal{G})$ and $LE^+_{CN}(\mathcal{G})$ analogous to the bound given in  \cite[Remark 3.8]{DM-2014}.

\begin{thm}\label{new bound} Let $\mathcal{G}$ be a graph of order $|V(\mathcal{G})|$ with $|e(\mathcal{G})|$ edges and the first Zagreb index $M_1(\mathcal{G})$. Then
$LE_{CN}(\mathcal{G})$ and $LE^+_{CN}(\mathcal{G})$ are bounded above by
\[
2\,\left (1-\frac{1}{|V(\mathcal{G})|} \right)\,\Big(M_1(\mathcal{G})-2\,|e(\mathcal{G})|\Big).
\]
\end{thm}

\begin{proof}
We have
\begin{align*}
&S_\alpha(\mathcal{G})\leq S_{|V(\mathcal{G})|}(\mathcal{G})=tr(\CNRS(\mathcal{G}))=M_1(\mathcal{G})-2\,|e(\mathcal{G})| \\
\text{ and }&\\
&S^+_\beta(\mathcal{G})\leq S^+_{|V(\mathcal{G})|}(\mathcal{G})=tr(\CNRS(\mathcal{G}))=M_1(\mathcal{G})-2\,|e(\mathcal{G})|.
\end{align*}
Therefore, by Lemma \ref{Finite-sum equality} with $1\leq \alpha\leq |V(\mathcal{G})|$, we obtain
\[
LE_{CN}(\mathcal{G})\leq 2\,\left (1-\frac{1}{|V(\mathcal{G})|} \right)\,\Big(M_1(\mathcal{G})-2\,|e(\mathcal{G})|\Big).
\]
Similarly, we get the bound for $LE^+_{CN}(\mathcal{G})$.
\end{proof}
Note that the bounds obtained in Theorem \ref{new bound} are better then the bounds obtained in \eqref{rem-eq-2}.  We conclude this section with another upper bounds for $LE_{CN}(\mathcal{G})$ and $LE^+_{CN}(\mathcal{G})$ analogous to the bound obtained in  \cite[Theorem 5.5]{DMG-2017}. The following lemma is useful in this regard.

\begin{lem}{\rm \cite[Lemma 5.1]{DMG-2017}} \label{E-sqrt-inequality}
Let $A$ be a real symmetric matrix of order $n$ and let $d_1,d_2,\dots,d_n$ be the diagonal entries of the matrix $A^2$. Then
$$E(A)\leq \sum\limits_{i=1}^n \sqrt{d_i}.$$
\end{lem}

\begin{thm}\label{CNRS-diagonality} Let $\mathcal{G}$ be a graph of order $|V(\mathcal{G})|$ with $|e(\mathcal{G})|$ edges and the first Zagreb index $M_1(\mathcal{G})$. Then
\begin{enumerate}
\item $LE_{CN}(\mathcal{G})\leq \sum\limits_{i=1}^{|V(\mathcal{G})|} \sqrt{\left (\CNRS(\mathcal{G})_{i,i}-\displaystyle{\frac{M_1(\mathcal{G})-2\,|e(\mathcal{G})|}{|V(\mathcal{G})|}}\right )^2+\sum\limits_{k=1,\,k\neq i}^{|V(\mathcal{G})|}\,|N_{\mathcal{G}}(v_i)\cap N_{\mathcal{G}}(v_k)|^2}$,
\item $LE^+_{CN}(\mathcal{G})\leq \sum\limits_{i=1}^{|V(\mathcal{G})|} \sqrt{\left (\CNRS(\mathcal{G})_{i,i}-\displaystyle{\frac{M_1(\mathcal{G})-2\,|e(\mathcal{G})|}{|V(\mathcal{G})|}}\right )^2+\sum\limits_{k=1,\,k\neq i}^{|V(\mathcal{G})|}\,|N_{\mathcal{G}}(v_i)\cap N_{\mathcal{G}}(v_k)|^2}$,
\end{enumerate}
where $\CNRS(\mathcal{G})_{i,i}=d_{\mathcal{G}}(v_i)\,\Big(m_{\mathcal{G}}(v_i)-1\Big)$ for $1 \leq i \leq |V(\mathcal{G})|$ are main diagonal elements of $\CNRS(\mathcal{G})$ such that $\CNRS(\mathcal{G})_{1,1} \geq \CNRS(\mathcal{G})_{2,2} \geq \cdots \geq \CNRS(\mathcal{G})_{|V(\mathcal{G})|,|V(\mathcal{G})|}$.
\end{thm}
\begin{proof}
(a) Let $M=\CNL(\mathcal{G})-\frac{tr(\CNRS(\mathcal{G}))}{|V(\mathcal{G})|}I_{|V(\mathcal{G})|}$, where $I_{|V(\mathcal{G})|}$ is the identity matrix of size $|V(\mathcal{G})|$. Then  $\spec(M)=\left\lbrace \nu_i-\frac{tr(\CNRS(\mathcal{G}))}{|V(\mathcal{G})|}~:~1\leq i\leq |V(\mathcal{G})| \right\rbrace$, where $\cnlspec(\mathcal{G})=\{\nu_1,\nu_2,\dots,\nu_{|V(\mathcal{G})|}\}$. We have
\begin{align*}
    M^2=&\left (\CNL(\mathcal{G})-\frac{tr(\CNRS(\mathcal{G}))}{|V(\mathcal{G})|}I_{|V(\mathcal{G})|}\right )^2\\
       =&(\CNL(\mathcal{G}))^2-\frac{2\,tr(\CNRS(\mathcal{G}))}{|V(\mathcal{G})|}\CNL(\mathcal{G})+\frac{(tr(\CNRS(\mathcal{G})))^2}{|V(\mathcal{G})|^2}I_{|V(\mathcal{G})|}.
\end{align*}
Therefore, the $i$-th diagonal element of $M^2$ is
\[
(M^2)_{i,i}=(\CNL(\mathcal{G}))^2_{i,i}-\frac{2\,tr(\CNRS(\mathcal{G}))}{|V(\mathcal{G})|}(\CNRS(\mathcal{G}))_{i,i}+\frac{(tr(\CNRS(\mathcal{G})))^2}{|V(\mathcal{G})|^2}.
\]
We have
\begin{align*}
    (\CNL(\mathcal{G}))^2=&\Big(\CNRS(\mathcal{G})-\CN(\mathcal{G})\Big)^2\\
                         =&\CNRS(\mathcal{G})^2-\CNRS(\mathcal{G})\CN(\mathcal{G})-\CN(\mathcal{G})\CNRS(\mathcal{G})+\CN(\mathcal{G})^2.
\end{align*}
Therefore,
\begin{align*}
    (\CNL(\mathcal{G}))_{i,i}^2=&(\CNRS(\mathcal{G}))^2_{i,i}+(\CN(\mathcal{G}))^2_{i,i}
                              =(\CNRS(\mathcal{G})_{i,i})^2+\sum\limits_{k=1,\,k\neq i}^{|V(\mathcal{G})|}(\CN(\mathcal{G})_{i,k})^2.
\end{align*}
Hence,
\[
(M^2)_{i,i} =  \left ( \CNRS(\mathcal{G})_{i,i}-\frac{tr(\CNRS(\mathcal{G}))}{|V(\mathcal{G})|} \right )^2+\sum\limits_{k=1,\,k\neq i}^{|V(\mathcal{G})|}(\CN(\mathcal{G})_{i,k})^2.
\]

\vspace*{3mm}

\noindent
Since $tr(\CNRS(\mathcal{G}))=M_1(\mathcal{G})-2\,|e(\mathcal{G})|$, by Lemma \ref{E-sqrt-inequality}, we obtain
\[
E(M)\leq \sum\limits_{i=1}^{|V(\mathcal{G})|} \sqrt{\left ( \CNRS(\mathcal{G})_{i,i}-\frac{M_1(\mathcal{G})-2\,|e(\mathcal{G})|}{|V(\mathcal{G})|} \right )^2+\sum\limits_{k=1,\,k\neq i}^{|V(\mathcal{G})|}\,|N_{\mathcal{G}}(v_i)\cap N_{\mathcal{G}}(v_k)|^2}.
\]
Hence the result follows noting that
$$LE_{CN}(\mathcal{G})=E(M).$$

(b) Let $N=\CNSL(\mathcal{G})-\frac{tr(\CNRS(\mathcal{G}))}{|V(\mathcal{G})|}I_{|V(\mathcal{G})|}$, where $I_{|V(\mathcal{G})|}$ is the identity matrix of size $|V(\mathcal{G})|$. Then  $\spec(N)=\left\lbrace \sigma_i-\frac{tr(\CNRS(\mathcal{G}))}{|V(\mathcal{G})|}~:~1\leq i\leq |V(\mathcal{G})| \right\rbrace$, where $\cnqspec(\mathcal{G})=\{\sigma_1,\sigma_2,\dots,\sigma_{|V(\mathcal{G})|}\}$. We have
\begin{align*}
    N^2=&\left (\CNSL(\mathcal{G})-\frac{tr(\CNRS(\mathcal{G}))}{|V(\mathcal{G})|}I_{|V(\mathcal{G})|}\right )^2\\
       =&(\CNSL(\mathcal{G}))^2-\frac{2tr(\CNRS(\mathcal{G}))}{|V(\mathcal{G})|}\CNSL(\mathcal{G})+\frac{(tr(\CNRS(\mathcal{G})))^2}{|V(\mathcal{G})|^2}I_{|V(\mathcal{G})|}.
\end{align*}
Therefore the $i$-th diagonal element of $N^2$ is
\[
(N^2)_{i,i}=(\CNSL(\mathcal{G}))^2_{i,i}-\frac{2tr(\CNRS(\mathcal{G}))}{|V(\mathcal{G})|}(\CNRS(\mathcal{G}))_{i,i}+\frac{(tr(\CNRS(\mathcal{G})))^2}{|V(\mathcal{G})|^2}.
\]
We have
\begin{align*}
    (\CNSL(\mathcal{G}))^2=&\Big(\CNRS(\mathcal{G})+\CN(\mathcal{G})\Big)^2\\
                         =&\CNRS(\mathcal{G})^2+\CNRS(\mathcal{G})\CN(\mathcal{G})+\CN(\mathcal{G})\CNRS(\mathcal{G})+\CN(\mathcal{G})^2.
\end{align*}
Therefore,
\[
    (\CNSL(\mathcal{G}))_{i,i}^2=\Big(\CNRS(\mathcal{G})_{i,i}\Big)^2+\sum\limits_{k=1}^{|V(\mathcal{G})|}(\CN(\mathcal{G})_{i,k})^2.
\]
Hence,
\begin{align*}
    (N^2)_{i,i}=&\left ( \CNRS(\mathcal{G})_{i,i}-\frac{tr(\CNRS(\mathcal{G}))}{|V(\mathcal{G})|} \right )^2+\sum\limits_{k=1}^{|V(\mathcal{G})|}(\CN(\mathcal{G})_{i,k})^2.
\end{align*}

\vspace*{3mm}

\noindent
Since $tr(\CNRS(\mathcal{G}))=M_1(\mathcal{G})-2\,|e(\mathcal{G})|$, by Lemma \ref{E-sqrt-inequality}, we have
\[
E(N)\leq \sum\limits_{i=1}^{|V(\mathcal{G})|} \sqrt{\left ( \CNRS(\mathcal{G})_{i,i}-\frac{M_1(\mathcal{G})-2\,|e(\mathcal{G})|}{|V(\mathcal{G})|} \right )^2+\sum\limits_{k=1}^{|V(\mathcal{G})|}\,|N_{\mathcal{G}}(v_i)\cap N_{\mathcal{G}}(v_k)|^2}.
\]
Hence, the result follows noting that
$$LE^+_{CN}(\mathcal{G})=E(N).$$
\end{proof}

\begin{thm} If $\mathcal{G}$ is a $r$-regular graph of order $|V(\mathcal{G})|$, then $LE_{CN}(\mathcal{G})$ and $LE^+_{CN}(\mathcal{G})$ are bounded by $$\sum\limits_{i=1}^{|V(\mathcal{G})|}\sqrt{\sum\limits_{k=1}^{|V(\mathcal{G})|}\,|N_{\mathcal{G}}(v_i)\cap N_{\mathcal{G}}(v_k)|^2}.$$
\end{thm}

\begin{proof} Since $\mathcal{G}$ is a regular graph, by Lemma \ref{1kin1}, we obtain
\begin{align*}
tr(\CNRS(\mathcal{G}))&=\sum\limits^{|V(\mathcal{G})|}_{i=1}\,\sum\limits^{|V(\mathcal{G})|}_{j=1,\,j\neq i}\,|N_{\mathcal{G}}(v_i)\cap N_{\mathcal{G}}(v_k)|\\[2mm]
&=\sum\limits^{|V(\mathcal{G})|}_{i=1}\,\,d_{\mathcal{G}}(v_i)\,\Big(m_{\mathcal{G}}(v_i)-1\Big)\\[2mm]
&=\sum\limits^{|V(\mathcal{G})|}_{i=1}\,r\,(r-1)=|V(\mathcal{G})|\,r\,(r-1)
\end{align*}
and
  $$\CNRS(\mathcal{G})_{i,i}=r\,(r-1)=\frac{tr(\CNRS(\mathcal{G}))}{|V(\mathcal{G})|},~~~~~1\leq i\leq |V(\mathcal{G})|.$$
From Theorem \ref{CNRS-diagonality}, we get the result.
\end{proof}


\section*{Acknowledgement} F. E. Jannat is supported by  DST INSPIRE Fellowship (IF200226). K. C. Das is supported by National Research Foundation funded by the Korean government (Grant No. 2021R1F1A1050646).

\vspace*{4mm}

\noindent
{\Large \bf Competing Interests:} The authors declare no conflict of interest.


\end{document}